\documentclass[11pt, a4paper]{article}
\usepackage{amsmath,amssymb,latexsym,graphics,epsfig}
\usepackage{hyperref}
\usepackage{color}
\usepackage{amsthm}

\newtheorem*{theo}{Theorem}
\newtheorem*{lema}{Lemma}
\newtheorem*{propo}{Proposition}

\newtheorem*{ogthm}{Odd-girth theorem}

\def\prova{{\boldmath  $Proof.$}\hskip 0.3truecm}

\def\final{\mbox{ \quad $\Box$}}

\def\A{\mbox{\boldmath $A$}}

\def\Re{\mathbb R}

\def\dist{\mathop{\rm dist }\nolimits}

\def\Ker{\mathop{\rm Ker }\nolimits}

\def\j{\mbox{\boldmath $j$}}

\def\A{\mbox{\boldmath $A$}}

\def\E{\mbox{\boldmath $E$}}
\def\Ei{{\cal E}}

\def\G{\Gamma}

\def\I{\mbox{\boldmath $I$}}

\def\J{\mbox{\boldmath $J$}}

\def\sp{\mbox{\rm sp}}

\def\Ker{\mathop{\rm Ker }\nolimits}

\def\sp{\mathop{\rm sp }\nolimits}

\begin{document}
\title{A Short Proof of the Odd-Girth Theorem
%\thanks{Research supported by the Ministerio de Educaci\'on y
%Ciencia, Spain, and the European Regional Development Fund under
%project MTM2008-06620-C03-01 and by the Catalan Research Council
%under project 2009SGR1387.}
}

\author{E.R. van Dam$^\ddag$ and M.A. Fiol$^\dag$
\\ \\
{\small $^\ddag$Tilburg University, Dept. Econometrics and O.R.} \\
{\small Tilburg, The Netherlands} {\small (e-mail: {\tt
edwin.vandam@uvt.nl})} \\
{\small $^\dag$Universitat Polit\`ecnica de Catalunya, Dept. de Matem\`atica Aplicada IV} \\
{\small Barcelona, Catalonia} {\small (e-mail: {\tt fiol@ma4.upc.edu})}}

\date{}

\maketitle

\begin{abstract}
\noindent Recently, it has been shown that a connected graph $\G$ with $d+1$
distinct eigenvalues and odd-girth $2d+1$ is distance-regular. The proof of
this result was based on the spectral excess theorem. In this note we present
an alternative and more direct proof which does not rely on the spectral excess
theorem, but on a known characterization of distance-regular graphs in terms of
the predistance polynomial of degree $d$.
\end{abstract}

%%%%%%%%%%%%%%%%%%%%%%%%%%%%%%%%%%%%%%%%%%%%%%%%
\section{Introduction}
%%%%%%%%%%%%%%%%%%%%%%%%%%%%%%%%%%%%%%%%%%%%%%%%
The spectral excess theorem \cite{fg97} states that a regular graph $\G$ is
distance-regular if and only if its spectral excess (a number which can be
computed from the spectrum of $\G$) equals its average excess (the mean of the
numbers of vertices at maximum distance from every vertex), see
\cite{vd08,fgg10} for short proofs. Using this theorem, Van Dam and Haemers
\cite{vdh11} proved the below odd-girth theorem for regular graphs.

\begin{ogthm} A connected graph with $d+1$ distinct
eigenvalues and finite odd-girth at least $2d+1$ is a distance-regular
generalized odd graph $($that is, a distance-regular graph with diameter $D$
and odd-girth $2D+1$\/$)$. \end{ogthm}

In the same paper, the authors posed the problem of deciding whether the
regularity condition is necessary or, equivalently, whether or not there are
nonregular graphs with $d+1$ distinct eigenvalues and odd girth $2d+1$.
Moreover, they proved this in the negative for the case $d+1=3$, and claimed to
have proofs for the cases $d+1\in \{4,5\}$. In a recent paper, Lee and Weng
\cite{lw11} used a variation of the spectral excess theorem for nonregular
graphs to show that, indeed, the regularity condition is not necessary. The
odd-girth theorem generalizes the result by Huang and Liu \cite{HL99} that
states that every graph with the same spectrum as a generalized odd graph must
be such a graph itself. Well-known examples of generalized odd graphs are the
odd graphs and the folded cubes.

In this note we give a short and direct proof of the more general result
without using any of the spectral excess theorems, but only a known
characterization of distance-regularity in terms of the predistance polynomial
$p_d$ of highest degree.

\section{Preliminaries}
Here we give some basic notation and results on which our proof of the
odd-girth theorem is based. For more background on spectra  of graphs,
distance-regular graphs, and their characterizations, see
\cite{b93,bcn89,bh12,cds82,dkt12,f02}.

Let $\G$ be a connected graph with vertex set $V$, order $n=|V|$, and adjacency
matrix $\A$. The spectrum of $\G$ (that is, of $\A$) is denoted by $\sp \G =
\{\lambda_0^{m_0},\lambda_1^{m_1},\dots, \lambda_d^{m_d}\}$, with distinct
eigenvalues $\lambda_0>\lambda_1>\cdots >\lambda_d$, and corresponding
multiplicities $m_i=m(\lambda_i)$. The {\em predistance polynomials} $p_i$
($i=0,1,\ldots,d$) of $\G$, form a sequence of orthogonal polynomials with
respect to the scalar product $\langle f,g\rangle=\frac{1}{n}\sum_{i=0}^d m_i
f(\lambda_i)g(\lambda_i)$, normalized in such a way that
$\|p_i\|^2=p_i(\lambda_0)$. Then, modulo the minimal polynomial of $\A$, these
polynomials satisfy a three-term recurrence relation of the form
\begin{equation}
\label{recur-pol}
xp_i=\beta_{i-1}p_{i-1}+\alpha_i p_i+\gamma_{i+1}p_{i+1}\qquad
(i=0,1,\dots,d),
\end{equation}
where we let $\beta_{-1}p_{-1}=0$ and $\gamma_{d+1}p_{d+1}=0$. Note that if
$\G$ is distance-regular and $\A_i$ stands for the distance-$i$ matrix, then
$\A_i=p_i(\A)$ for $i=0,1,\ldots,d$. Our proof of the odd-girth theorem relies
mainly on the following result which was first proved in \cite{fgy96b} (see
also \cite{vd08}, \cite{fgg10}).

\begin{propo}
\label{propo-charac-drg} A regular graph $\G$ with $d+1$ distinct eigenvalues
is distance-regular if and only if the predistance polynomial $p_d$ satisfies
$p_d(\A)=\A_d$.
\end{propo}

We remark that it is fairly easy to prove this characterization by backward
induction, using the recurrence relation \eqref{recur-pol}, the fact that the
Hoffman polynomial $H$ equals $p_0+p_1+\cdots+p_d$, and that for regular graphs
$H(\A)$ equals the all-1 matrix $\J$.

By the matrices  $\E_i$ we denote the {\it $($principal\/$)$ idempotents} of
$\A$ representing the orthogonal projections of $\Re^n$ onto the eigenspaces
$\Ei_i=\Ker (\A-\lambda_i \I)$, for $i=0,1,\dots,d$. In particular, if $\G$ is
regular, then the all-$1$ vector $\j$ is a $\lambda_0$-eigenvector and
$\E_0=\frac{1}{n}\j\j^{\top}=\frac{1}{n}\J$. The diagonal entries of these
idempotents, $m_{u}(\lambda_i)=(\E_i)_{uu}$, have been called the {\em
$u$-local multiplicities} of the eigenvalue $\lambda_i$. Graphs for which these
local multiplicities are independent of the vertex $u$ (that is, for which
every idempotent has constant diagonal) are called {\em spectrum-regular}. The
local multiplicities allow us to compute the number of closed $\ell$-walks from
$u$ to itself in the following way:
\begin{equation}
\label{crossed-mul->num-walks} a_{u}^{({\ell})} =(\A^{\ell})_{uu}=
\sum_{i=0}^d m_{u}(\lambda_i)\lambda_i^{\ell} \qquad (\ell =0,1,2,\ldots).
\end{equation}

A graph is {\em walk-regular} (a concept introduced by Godsil and McKay
\cite{gmk80}) if the number $a_{u}^{({\ell})}$ of closed walks of length $\ell$
does not depend on $u$, for every $\ell=0,1,2,\ldots$. Clearly, a graph is
walk-regular if and only if it is spectrum-regular, and every walk-regular
graph is regular; properties that will be used in our proof of the odd-girth
theorem.

\section{The proof}

Now let us consider a connected graph $\G$ with $d+1$ distinct eigenvalues and
finite odd-girth (at least) $2d+1$, and the corresponding predistance
polynomials with recurrence \eqref{recur-pol}. As was shown in \cite{vdh11} by
an easy inductive argument, in this particular case we have that $\alpha_i=0$
for $i=0,1,\ldots,d-1$ and the polynomials $p_i$ are even or odd depending on
$i$ being even or odd, respectively. Moreover, $\alpha_d \neq 0$ (even though
Van Dam and Haemers \cite{vdh11} restrict to regular graphs, the regularity
condition is not used by them; the argument is also implicitly used by Lee and
Weng \cite{lw11}). In order to prove the odd-girth theorem, we first need the
following lemma.

\begin{lema}
\label{plusminuslemma}
Let $\G$ be a connected graph with $d+1$ distinct
eigenvalues and odd-girth $2d+1$. If $\lambda$ is an eigenvalue of $\G$, then
$-\lambda$ is not. In particular, all eigenvalues are nonzero.
\end{lema}
\prova Assume that both $\lambda$ and $-\lambda$ are eigenvalues of $\G$, that
is, they are both roots of the minimal polynomial. That means that we can plug
in $\lambda$ and $-\lambda$ in the recurrence relations \eqref{recur-pol} and,
in particular, we obtain the two equations $\pm \lambda p_d(\pm
\lambda)=\beta_{d-1}p_{d-1}(\pm \lambda)+\alpha_dp_d(\pm \lambda)$. By using
that the predistance polynomials are odd or even as indexed, and that $\alpha_d
\neq 0$, it follows that $p_d(\lambda)=0$ (also in the case that $\lambda=0$),
which is a contradiction (because by the recurrence relations this would imply
that $p_i(\lambda)=0$ for all $i$, including $i=0$, but $p_0=1$).  \final

Now we are ready to prove the general setting of the odd-girth theorem without
using the spectral excess theorem.

\begin{theo}
\label{main-theo}
A connected graph $\G$ with $d+1$ distinct eigenvalues and odd-girth $2d+1$ is distance-regular.
\end{theo}
\prova Fist, let us prove that $\G$ is spectrum-regular (or walk-regular).
Since the number of odd cycles with length at most $2d-1$ is zero we have,
using \eqref{crossed-mul->num-walks},
\begin{equation*}
\sum_{i=1}^d m_u(\lambda_i)\lambda_i^{2\ell-1}=-m_u(\lambda_0)\lambda_0^{2\ell-1} \qquad (\ell=1,2,\ldots,d).
\end{equation*}
This can be seen as a determined system of $d$ equations and $d$ unknowns
$m_u(\lambda_i)$ ($i=1,2,\ldots,d$). Indeed, by the properties of Vandermonde
matrices and the above lemma, the determinant of its coefficient matrix is
\begin{eqnarray*}
\left|
\begin{array}{cccc}
\lambda_1  & \lambda_2 & \cdots & \lambda_d \\
\lambda_1^3  & \lambda_2^3 & \cdots & \lambda_d^3 \\
\vdots &\vdots & \ddots &  \vdots \\
\lambda_1^{2d-1}  & \lambda_2^{2d-1} &\cdots & \lambda_d^{2d-1}
\end{array}
\right| & =& \prod_{i=1}^d\lambda_i\left|\begin{array}{cccc}
1  & 1 & \cdots & 1 \\
\lambda_1^2  & \lambda_2^2 & \cdots & \lambda_d^2 \\
\vdots &\vdots & \ddots &  \vdots \\
\lambda_1^{2d-2}  & \lambda_2^{2d-2} &\cdots & \lambda_d^{2d-2}
\end{array}\right| \\
 & = & \prod_{i=1}^d\lambda_i\prod_{d\ge i>j\ge 1}(\lambda_i^2-\lambda_j^2) \neq 0.
\end{eqnarray*}
Thus, there exist constants $\alpha_i$ such that $m_u(\lambda_i)=\alpha_i
m_u(\lambda_0)$, for $i=0,1,\ldots,d$. From this it follows that
$m_u(\lambda_0)\sum_{j=0}^d \alpha_j=\sum_{j=0}^dm_u(\lambda_j)=1$, where the
last equality follows from the fact that the sum of all idempotents equals the
identity matrix. Thus, for every $i=0,1,\ldots,d$,
$m_u(\lambda_i)=\alpha_i/\sum_{j=0}^d \alpha_j$, which does not depend on $u$,
and hence $\G$ is spectrum-regular (and walk-regular).

Next, let us show that $p_d(\A)=\A_d$. Since $\G$ is regular, the Hoffman
polynomial $H=p_0+p_1+\cdots+p_d$ satisfies $H(\A)=\J$ and hence
$(p_d(\A))_{uv}=1$ if $\dist(u,v)=d$. Besides, from the parity of the
predistance polynomials, it follows that $(p_i(\A))_{uv}=0$ if $\dist(u,v)$ and
$i$ have different parity (otherwise, $\G$ would have an odd cycle of length
smaller than $2d+1$). So $(p_d(\A))_{uv}=0$ for every pair of vertices $u,v$
whose distance has a different parity than $d$. If $\dist(u,v)$ is smaller than
$d$, but with the same parity, then from the recurrence \eqref{recur-pol} we
get
$$
(\A p_d(\A))_{uv}=\beta_{d-1}(p_{d-1}(\A))_{uv}+\alpha_d (p_d(\A))_{uv}=
\alpha_d (p_d(\A))_{uv}
$$
(because $\dist(u,v)$ and $d-1$ have different parity). But the first term is
$$
\sum_{w\in V} (\A)_{uw} (p_d(\A))_{wv}=\sum_{w\in \G(u)}(p_d(\A))_{wv}=0
$$
since $\dist(w,v)=\dist(u,v)\pm 1$ has a different parity than $d$. Thus, as
$\alpha_d\neq 0$, we find that also in this case $(p_d(\A))_{uv}=0$.
Consequently, $p_d(\A)=\A_d$ and by the above proposition, $\G$ is
distance-regular. \final

%%%%%%%%%%%%%%%%%%%%%%%%%%%%%%%%%%%%%%%%%%%%%%%%
%Bibliografia
%%%%%%%%%%%%%%%%%%%%%%%%%%%%%%%%%%%%%%%%%%%%%%%%

\end{document}